\documentclass[10pt]{amsart}
\usepackage{amssymb}
\usepackage{latexsym}
\usepackage{enumerate}

\newtheorem{thm}{Theorem}[section]
\newtheorem{defin}[thm]{Definition}
\newtheorem{cor}[thm]{Corollary}
\newtheorem{prop}[thm]{Proposition}
\newtheorem{lemma}[thm]{Lemma}
\newtheorem{eg}[thm]{Example}
\newtheorem{rmk}[thm]{Remark}

\newtheorem{conjecture}[thm]{Conjecture}

\newcommand{\integers}{\ensuremath{\mathbb{Z}}}
\newcommand{\rationals}{\ensuremath{\mathbb{Q}}}
\newcommand{\reals}{\ensuremath{\mathbb{R}}}
\newcommand{\complex}{\ensuremath{\mathbb{C}}}

\newcommand{\finfield}[1]{\ensuremath{\mathbb{F}_{#1}}}

\newcommand{\pf}[1][]{\emph{Proof#1:\ }}

\newcommand{\abs}[1]{\ensuremath{\left|#1\right|}}
\newcommand{\defsetshort}[1]{\ensuremath{\left\{#1\right\}}}

\newcommand{\defsetspan}[1]{\ensuremath{<\!\!#1\!\!>}}

\renewcommand{\mod}[1]{\ensuremath{\ \left(\mbox{mod }#1\right)}}

\newcommand{\symgroup}[1]{\ensuremath{\mathbb{S}_{#1}}}
\newcommand{\altgroup}[1]{\ensuremath{\mathbb{A}_{#1}}}
\newcommand{\dihgroup}[1]{\ensuremath{\mathbb{D}_{#1}}}
\newcommand{\cyclgroup}[1]{\ensuremath{\integers_{#1}}}

\newcommand{\bdihgroup}[1]{\ensuremath{\overline{\mathbb{D}}_{#1}}}

\newcommand{\glgroup}[2]{\ensuremath{\mbox{GL}_{#1}\left(#2\right)}}
\newcommand{\slgroup}[3][]{\ensuremath{\mbox{SL}^{#1}_{#2}\left(#3\right)}}
\newcommand{\pglgroup}[2]{\ensuremath{\mbox{PGL}_{#1}\left(#2\right)}}
\newcommand{\pslgroup}[2]{\ensuremath{\mbox{PSL}_{#1}\left(#2\right)}}

\newcommand{\spgroup}[2]{\ensuremath{\mbox{Sp}_{#1}\left(#2\right)}}
\newcommand{\proj}[2][]{\ensuremath{\mathbb{P}_{#1}^{#2}}}

\newcommand{\deffunname}[3]{\ensuremath{#1:#2\rightarrow#3}}

\newcommand{\defmap}[2]{\ensuremath{#1\mapsto#2}}

\newcommand{\canclass}[1]{\ensuremath{\mbox{K}_{#1}}}

\title{Five-dimensional weakly-exceptional quotient singularities.}
\author{Dmitrijs Sakovics}

\begin{document}
\begin{abstract}
A singularity is said to be weakly-exceptional if it has a unique purely log terminal blow up. This is a natural generalization of the surface singularities of types $D_{n}$, $E_{6}$, $E_{7}$ and $E_{8}$.  Since this idea was introduced, quotient singularities of this type have been classified in dimensions up to most $4$. This note extends that classification to dimension $5$.
\end{abstract}
\maketitle

\section{Introduction}
Let $G\subset\glgroup{n}{\complex}$ be a finite group. It makes sense to study the quotient singularities on the varieties of the form $\complex^{n}/G$ (from now on, these will be referred to as the singularities induced by $G$). When studying singularities (and, in particular, quotient singularities), one may consider the following type of birational morphisms:
\begin{thm}[see~{\cite[Theorem~3.7]{Cheltsov-Shramov09}}]\label{thm:plt}
Let $(V\ni O)$ be a~germ of a~Kawamata log terminal singularity. Then there
exists a birational morphism \deffunname{\pi}{W}{V} such that the following
hypotheses are satisfied:
\begin{itemize}
\item the~exceptional locus of $\pi$ consists of one irreducible divisor
$E$ such that $O\in\pi(E)$,%
\item the~log pair $(W,E)$ has purely log terminal singularities.
\item the~divisor $-E$ is a~$\pi$-ample \rationals-Cartier divisor.
\end{itemize}
\end{thm}
\begin{defin}[{\cite{kudryavtsev01}}]
Let $(V\ni O)$ be a~germ of a~Kawamata log terminal singularity,
and \deffunname{\pi}{W}{V} be~a~birational morphism satisfying the conditions
of Theorem~\ref{thm:plt}. Then $\pi$ is a~\emph{plt blow-up} of
the~singularity.
\end{defin}
This naturally leads to the following definition:
\begin{defin}[{\cite{kudryavtsev01}}]
We say that the singularity $(V\ni O)$ is \emph{weakly-exceptional} if it has a unique plt blow-up.
\end{defin}
\begin{rmk}\label{rmk:sl5:dim2}
 This definition naturally generalizes the properties of quotients of $\complex^{2}$ by the action of binary dihedral (also known as dicyclic), tetrahedral, octahedral and icosahedral groups into higher dimensions.
\end{rmk}

One has the following criterion for a singularity being weakly-exceptional:
\begin{thm}[see {\cite[Theorem~3.15]{Cheltsov-Shramov09}}]\label{thm:geom_criterion}
Take $G\subset\glgroup{N}{\complex}$ with no quasi-reflections, and let $\bar{G}$ be its natural projection into \pglgroup{N}{\complex}. Then the singularity $\complex^{N}/G$ is weakly-exceptional if and only if the pair $\left(\proj{N-1},\Delta\right)$ is log canonical for any $\bar{G}$-invariant effective \rationals-divisor $\Delta\sim_{\rationals}-\canclass{\proj{N-1}}$.
\end{thm}
\begin{cor}\label{cor:sl5:low_invars}
 If $G\subset\glgroup{N}{\complex}$ with no quasi-reflections has a semi-invariant of degree at most $N-1$, then the singularity induced by it is not weakly-exceptional.
\end{cor}
\begin{rmk}
 The reverse implication does not hold in general, for example for $N=4$ (see~\cite{Sakovics10}).
\end{rmk}
\begin{thm}[{\cite[Theorem~1.30]{Cheltsov-Shramov09}}]\label{thm:WExc-trans}
Let $G\subset\glgroup{N}{\complex}$ be a finite subgroup containing no quasi-reflections that induces a
weakly-exceptional singularity. Then $G$ is irreducible.
\end{thm}

In fact, in dimension $2$, the induced quotient singularity is weakly-exceptional exactly when the group action is irreducible. Unfortunately, this fails already in dimension $3$.

It follows from the Chevalley--Shephard--Todd Theorem (see~\cite[Theorem~4.2.5]{Springer77}) that to study the weak exceptionality of $\complex^{n}/G$, one can always assume that $G$ contains no quasi-reflections. Moreover, it follows from Theorem~\ref{thm:geom_criterion} that the weak exceptionality only depends on the image of $G$ under the natural projection to \pglgroup{n}{\complex}. So to study the weak exceptionality of $\complex^{n}/G$, it is enough to consider the case of  $G\subset\slgroup{n}{\complex}$.

The classification of the groups giving rise to weakly-exceptional singularities in dimension $2$ is well-known:
\begin{thm}[rephrasing~{\cite[Section~5.2.3]{Shokurov92-Eng}}]
 Let $G\subset\slgroup{2}{\complex}$ be a finite group. Then $G$ induces a weakly-exceptional singularity if and only if it is a non-abelian binary dihedral, tetrahedral, octahedral or icosahedral group.
\end{thm}
The groups giving rise to weakly-exceptional singularities in dimensions $3$ and $4$ have recently been classified (see~\cite{Sakovics10}). Due to the large number of irreducible groups in dimensions higher than $2$, it makes more sense to look at the irreducible groups that give rise to non-weakly-exceptional singularities. In particular, in dimension $3$ only finitely many conjugacy classes do so. Unfortunately, the same is not true in dimension $4$:
\begin{eg}\label{eg:dim4-infinite}
Write the coordinates of $\complex^{4}$ as a $2\times2$ matrix, and act on it by left and right multiplication by the elements of binary dihedral groups $\bdihgroup{2k},\bdihgroup{2l}\subset\slgroup{2}{\complex}$. Thus one gets an irreducible action of an arbitrarily large finite group, that has a semi-invariant quadric defined by the determinant of the matrix. This in turn implies (by Corollary~\ref{cor:sl5:low_invars}) that the induced singularity is not weakly-exceptional. For details, see~\cite{Sakovics10}.
\end{eg}
The purpose of this note is to prove the following result:
\begin{thm}\label{thm:sl5:fin-size}
 If $G\subset\slgroup{5}{\complex}$ is an irreducible monomial group that induces a non-weakly-exceptional singularity, then $\abs{G}\leq5\cdot4^{4}\cdot5!$, with this bound attained.
\end{thm}
Keeping in mind Remark~\ref{rmk:sl5:dim2} and the results of~\cite{Sakovics10}, one gets the following corollaries:
\begin{cor}
Take $p\in\defsetshort{2,3,5}$. Suppose $G\subset\slgroup{p}{\complex}$ is a finite subgroup acting irreducibly and monomially, but the singularity induced by $G$ is not weakly-exceptional. Then $\abs{G}\leq p\cdot\left(p-1\right)^{\left(p-1\right)}\cdot p!$.
\end{cor}
\begin{cor}
 Take $p\in\defsetshort{2,3,5}$. There are only finitely many finite groups $G\subset\slgroup{p}{\complex}$ (up to conjugation), such that $G$ acts irreducibly, but the singularity induced by $G$ is not weakly-exceptional.
\end{cor}

The proof of Theorem~\ref{thm:sl5:fin-size} relies on this note's main technical result, which one can consider to be the structure theorem for the irreducible groups in \slgroup{5}{\complex} inducing non-weakly-exceptional singularities (using notation introduced in Definition~\ref{defin:groupTypes} throughout):
\begin{thm}\label{thm:sl5:main}
 Let $G\subset\slgroup{5}{\complex}$ be a finite subgroup acting irreducibly. Then the singularity of $\complex^{5}/G$ is weakly-exceptional exactly when:
\begin{enumerate}
 \item\label{thm:sl5:main:prim} The action of $G$ is primitive and $G$ contains a subgroup isomorphic to the Heisenberg group of all unipotent $3\times3$ matrices over \finfield{5} (for a better classification of all such groups, see~\cite{Horrocks-Mumford73}).
 \item\label{thm:sl5:main:mono} The action of $G$ is monomial (making $G\cong D\rtimes T$, with $D$ an abelian group as above and $T$ a transitive subgroup of \symgroup{5}), and (using notation from Section~\ref{sect:sl5:mono}) none of the following hold:
\begin{itemize}
 \item $D$ is central in \slgroup{5}{\complex}. In this case, $G$ can be isomorphic to \altgroup{5}, \symgroup{5}, or their central extensions by \cyclgroup{5}. 
 \item $\abs{G}=55$ or $55\cdot5$ with $\abs{D}=11$ or $11\cdot5$ resp., $T\cong\cyclgroup{5}\subset\symgroup{5}$, and there is a $k\in\integers$, $1\leq k\leq4$, such that $D$ is generated by $\left[11,1,4^{k},4^{2k},4^{3k},4^{4k}\right]$ and (in the latter case) also the scalar element $\zeta_{5}\!\cdot\!\mbox{Id}$. In this case, $G$ is isomorphic to $\cyclgroup{11}\rtimes\cyclgroup{5}$ or $\left(\cyclgroup{5}\times\cyclgroup{11}\right)\rtimes\cyclgroup{5}$.
 \item $\abs{G}=305$ or $305\cdot5$ with $\abs{D}=61$ or $61\cdot5$ resp., $T\cong\cyclgroup{5}\subset\symgroup{5}$, and there is a $k\in\integers$, $1\leq k\leq4$, such that $D$ is generated by $\left[61,1,34^{k},34^{2k},34^{3k},34^{4k}\right]$ and (in the latter case) also the scalar element $\zeta_{5}\!\cdot\!\mbox{Id}$. In this case, $G$ is isomorphic to $\cyclgroup{61}\rtimes\cyclgroup{5}$ or $\left(\cyclgroup{5}\times\cyclgroup{61}\right)\rtimes\cyclgroup{5}$.
 \item There exists some $d\in\defsetshort{2,3,4}$ and $\omega$ with $\omega^{5}=1$, such that:
 \begin{itemize}
  \item $\forall g\in D$, $g^{d}$ is a scalar.
  \item $\abs{D}\in\defsetshort{d^{k},5\cdot d^{k}}$ (depending on whether $D$ contains any non-trivial scalar elements) with $1\leq k\leq4$.
  \item The polynomial $x_{1}^{d}+\omega x_{2}^{d}+\omega^{2}x_{3}^{d}+\omega^{3}x_{4}^{d}+\omega^{4}x_{5}^{d}$ is $G$-semi-invariant.
 \end{itemize}
\end{itemize}
\end{enumerate}
\end{thm}
\pf Let $G\subset\slgroup{5}{\complex}$ be a finite group. Since $5$ is a prime, $G$ is either primitive or monomial. This means that the result follows immediately from Lemma~\ref{lemma:sl5:prim} and the considerations in Section~\ref{sect:sl5:mono}.\qed

\pf[{ of Theorem~\ref{thm:sl5:fin-size}}] Follows directly from Theorem~\ref{thm:sl5:main}. The bound is attained by a group $G=D\rtimes T$ with $D=\cyclgroup{5}\times\cyclgroup{4}^{4}$ acting by scalar multiplication of coordinates of $\complex^{5}$, and $T\cong\symgroup{5}$ acting by permuting the basis. Here, $\mbox{Z}\left(G\right)=\cyclgroup{5}$. This group preserves the polynomial $\sum_{i=1}^{5}{x_{i}^{4}}$.\qed

This leads to the following conjecture:
\begin{conjecture}\label{conj:sl5:fin-number}
 For any prime $p$, there are only finitely many finite groups $G\subset\slgroup{p}{\complex}$ (up to conjugation), such that $G$ acts irreducibly, but the singularity induced by $G$ is not weakly-exceptional.
\end{conjecture}
It seems that an even stronger result holds: take any prime $p$ and suppose $G\subset\slgroup{p}{\complex}$ is a finite subgroup acting irreducibly and monomially, but the singularity induced by $G$ is not weakly-exceptional. Then $\abs{G}\leq p\cdot\left(p-1\right)^{\left(p-1\right)}\cdot p!$.

Note that Conjecture~\ref{conj:sl5:fin-number} can easily be shown to fail for infinitely many composite dimensions, as the construction in Example~\ref{eg:dim4-infinite} can easily be generalised to any dimension $n=k^{2}$.

\section{General considerations}
\begin{defin}[{\cite[\S2]{Hofling}}]\label{defin:groupTypes}
 Given a representation of a group $G$ on a space $V$, a \emph{system of imprimitivity} for the action is a set \defsetshort{V_{1},\ldots,V_{k}} of distinct subspaces of $V=V_{1}\oplus\cdots\oplus V_{k}$, such that $\forall i$ and $\forall g\in G$, $\exists j$ with $g\left(V_{i}\right)=V_{j}$. Clearly, \defsetshort{V} will always be one such system. If this is the only system of imprimitivity for this action, this action is called \emph{primitive}. If there is a system, where all the $V_{i}$-s are $1$-dimensional, then the action is called \emph{monomial}. If for any system of imprimitivity \defsetshort{V_{1},\ldots,V_{k}} and any $1\leq i,j\leq k$, $\exists g_{i,j}\in G$, such that $g_{i,j}\left(V_{i}\right)=V_{j}$, then the action is called irreducible.
\end{defin}

Since any group $G\subset\glgroup{5}{\complex}$ comes with a canonical faithful representation, it makes sense to say that the group itself, rather than that representation is primitive, monomial or irreducible.

\begin{thm}[\cite{Cheltsov-Shramov11}]\label{thm:classification}
Let $G$ be a finite subgroup in \glgroup{5}{\complex} that does not contain reflections. Then the~singularity $\complex^5/G$ is weakly-exceptional if and only if the~group $G$ is irreducible and does not have semi-invariants of degree at most $4$.
\end{thm}
It is worth noting that the property only depends on the projection of $G$ into \pglgroup{5}{\complex}. Therefore, from now on it will be assumed that $G\subset\slgroup{5}{\complex}$. If it is not, take instead a group $G'\subset\slgroup{5}{\complex}$ that has the same projection into \pglgroup{5}{\complex}.

This theorem provides two possible approaches to computing the list of irreducible groups giving rise to singularities in dimension $5$ that are not weakly-exceptional: either by obtaining a list of finite groups of automorphisms of projective threefolds of low degrees and seeing which of their actions descend to actions on \proj{4}, or by directly computing which groups have semi-invariant polynomials of degree at most $4$ in $5$ variables. Since the finite subgroups of \slgroup{5}{\complex} fit into two small families, that are relatively easy to work with, it has been chosen to follow the second approach.

To begin, it is easier to deal with the case of $G$ being a primitive group first, and then look into the monomial case.
\begin{thm}[{\cite[\S8.5]{Feit70}}]\label{thm:sl5:primList}
 If $G\subset\slgroup{5}{\complex}$ is a finite group acting primitively, then either $G$ is one of \altgroup{5},\altgroup{6}, \symgroup{5}, \symgroup{6}, \pslgroup{2}{11} and \spgroup{4}{\finfield{3}},
 or $G$ is a subgroup of the normalizer $\mathbb{HM}$ of the Heisenberg group $\mathbb{H}$ of all unipotent $3\times3$ matrices over \finfield{5}, such that $\mathbb{H}\subset G\subseteq\mathbb{HM}$.
\end{thm}


\begin{lemma}\label{lemma:sl5:prim}
Let $G\subset\slgroup{5}{\complex}$ be a finite primitive subgroup. Then $G$ gives rise to a weakly-exceptional singularity if and only if it contains a subgroup isomorphic to the Heisenberg group $\mathbb{H}$.
\end{lemma}
\pf Since there is a very small number of such groups (see Theorem~\ref{thm:sl5:primList}), one can simply look at the low symmetric powers of their $5$-dimensional irreducible representations. This gives:
\begin{itemize}
 \item The actions of \altgroup{5}, \symgroup{5}, \altgroup{6}, \symgroup{6} have semi-invariants of degree $2$, since they are conjugate to sub-
groups of \glgroup{5}{\reals}
 \item The action of \pslgroup{2}{11} has a semi-invariant of degree $3$, the Klein cubic threefold (see~\cite{Adler78}).
 \item The action of \spgroup{4}{\finfield{3}} has a semi-invariant of degree $4$, the Burkhardt quartic threefold (see~\cite{Burkhardt1891}).
 \item If $G$ contains the Heisenberg group $\mathbb{H}$, then $G$ cannot have any semi-invariants of degree at most $4$ (either apply Theorem~\ref{thm:classification} to~\cite[Theorem~1.15]{Cheltsov-Shramov11} or apply Lemma~\ref{lemma:sl5:elts_of_D} to the (monomial) representations of $H$ of dimension at most $5$).
\end{itemize}

\section{Monomial groups}\label{sect:sl5:mono}
Throughout this section, $\zeta_{n}$ will be used to denote a primitive $n$-th root of unity. This will be chosen consistently for different $n$, i.e.\ so that $\zeta_{mn}^{m}=\zeta_{n}$.

If $G\subset\slgroup{5}{\complex}$ is an finite irreducible monomial group, then take its system of imprimitivity consisting of $1$-dimensional subspaces. Let $D$ be the normal subgroup of $G$ preserving these subspaces. Then clearly, $D$ is abelian, and $G=D\rtimes T$, where $T$ is a transitive subgroup of \symgroup{5} permuting the spaces. Moreover, there is a basis for $\complex^{5}$ in which $D$ acts by multiplication by diagonal matrices, and there is an element $\tau\in G\setminus D$ acting by \defmap{\left(x_{1},x_{2},x_{3},x_{4},x_{5}\right)}{\left(x_{2},x_{3},x_{4},x_{5},x_{1}\right)}.

To establish non-ambiguous notation, one needs to mention that in this paper the notation \dihgroup{2n} will mean the dihedral group of $2n$ elements, and $\mathbb{GA}\left(1,5\right)=\defsetspan{\left(1\ 2\ 3\ 4\ 5\right),\left(2\ 3\ 5\ 4\right)}\subset\symgroup{5}$ is the General Affine group with parameters $\left(1,5\right)$. Furthermore, for any $g\in G$ and any polynomial $f$ write $g\left(f\right)=f\circ g$.

\begin{rmk}[{See, for example, appendix of~\cite{Swallow}}]\label{rmk:sl5:subgroups_of_S5}
 If $G$ is not generated by $D$ and $\tau$, then $\cyclgroup{5}\subsetneq T\subseteq\symgroup{5}$, so it is a well-known fact that $T$ must be one of \dihgroup{10}, $\mathbb{GA}\left(1,5\right)$, \altgroup{5} and \symgroup{5}, (up to choosing $\tau$) generated by $\left(1\ 2\ 3\ 4\ 5\right)$ (corresp.\ to $\tau$) and $\left(2\ 5\right)\left(3\ 4\right)$, $\left(2\ 3\ 5\ 4\right)$, $\left(1\ 2\ 3\right)$ or $\left(1\ 2\right)$ respectively.
\end{rmk}



Since $G$ is a finite group, any $g\in D$ must be multiplying the coordinates by roots of unity. From now on, write $\left[n,a_{1},a_{2},a_{3},a_{4},a_{5}\right]$ for the element acting as \defmap{\left(x_{1},x_{2},x_{3},x_{4},x_{5}\right)}{\left(\zeta_{n}^{a_{1}}x_{1},\zeta_{n}^{a_{2}}x_{2},\zeta_{n}^{a_{3}}x_{3},\zeta_{n}^{a_{4}}x_{4},\zeta_{n}^{a_{5}}x_{5}\right)}.
It is clear that $\left[n,a_{1},a_{2},a_{3},a_{4},a_{5}\right]=\left[kn,ka_{1},ka_{2},ka_{3},ka_{4},ka_{5}\right]$ for any $k\in\integers_{>0}$, so it will always be assumed that the presentation has the minimal possible $n\in\integers_{>0}$. Note that since $g\in\slgroup{5}{\complex}$, it must be true that $\sum_{i}{a_{i}}=nk$ for some $k\in\integers$. Also note that replacing $a_{i}$ by $a_{i}\pm n$ gives the same element.
\begin{lemma}\label{lemma:sl5:allscalars}
 If all the elements of $D$ are scalar, then the singularity induced by $G$ is not weakly-exceptional.
\end{lemma}
\pf In this case, $G$ must be either one of the groups mentioned in Remark~\ref{rmk:sl5:subgroups_of_S5} or a central extension of one of them by \cyclgroup{5}. On this list, the only groups that have irreducible $5$-dimensional representations are $\altgroup{5}$, $\symgroup{5}$ and their central extensions by \cyclgroup{5}. It is easy to see that all of these have semi-invariants of degree $2$.\qed

From now on, one can assume that $D$ contains a non-scalar element.

\begin{lemma}\label{lemma:sl5:makeprime}
 Let $g\in D$ be a non-scalar element of order $pq$ for some integers $p,q>1$. Then either $p=5$, or $\exists g'\in D$ a non-scalar element of order $p$.
\end{lemma}
\pf Set $g'=g^{q}$. Scalar elements in \slgroup{5}{\complex} have orders $1$ or $5$, so either $p=5$ or $g'$ is not a scalar.\qed

\begin{prop}
 Define the following monomials in $5$ variables $x_{1},\ldots,x_{5}$:
\[
\begin{array}{llll}
m_{1,1}=x_{1}&&&\\
\hline
m_{2,1}=x_{1}^{2}&
m_{2,2}=x_{1}x_{2}&
m_{2,3}=x_{1}x_{3}&\\
\hline
m_{3,1}=x_{1}^{3}&
m_{3,2}=x_{1}^{2}x_{2}&
m_{3,3}=x_{1}^{2}x_{3}&
m_{3,4}=x_{1}^{2}x_{4}\\
m_{3,5}=x_{1}^{2}x_{5}&
m_{3,6}=x_{1}x_{2}x_{3}&
m_{3,7}=x_{1}x_{2}x_{4}&\\
\hline
m_{4,1}=x_{1}^{4}&
m_{4,2}=x_{1}^{3}x_{2}&
m_{4,3}=x_{1}^{3}x_{3}&
m_{4,4}=x_{1}^{3}x_{4}\\
m_{4,5}=x_{1}^{3}x_{5}&
m_{4,6}=x_{1}^{2}x_{2}^{2}&
m_{4,7}=x_{1}^{2}x_{3}^{2}&
m_{4,8}=x_{1}^{2}x_{2}x_{3}\\
m_{4,9}=x_{1}^{2}x_{2}x_{4}&
m_{4,10}=x_{1}^{2}x_{2}x_{5}&
m_{4,11}=x_{1}^{2}x_{3}x_{4}&
m_{4,12}=x_{1}^{2}x_{3}x_{5}\\
m_{4,13}=x_{1}^{2}x_{4}x_{5}&
m_{4,14}=x_{1}x_{2}x_{3}x_{4}&&\\
\end{array}
\]
Then any polynomial $f$ of degree at most $4$ that is semi-invariant under the action of $\tau$ must be one of
\[\begin{array}{ll}
 A_{1}\sum_{j=0}^{4}{\omega^{j}\tau^{j}\left(m_{1,1}\right)} & \sum_{i=1}^{3}{\left[B_{i}\sum_{j=0}^{4}{\omega^{j}\tau^{j}\left(m_{2,i}\right)}\right]} \\
 \sum_{i=1}^{7}{\left[C_{i}\sum_{j=0}^{4}{\omega^{j}\tau^{j}\left(m_{3,i}\right)}\right]} & \sum_{i=1}^{14}{\left[D_{i}\sum_{j=0}^{4}{\omega^{j}\tau^{j}\left(m_{4,i}\right)}\right]}
\end{array}\]
where $A_{1},B_{i},C_{i},D_{i}\in\complex$ and $\omega$ is some (not necessarily primitive) fifth root of $1$.
\end{prop}
\pf The polynomial $f$ is semi-invariant under the action of $\tau$, so set $\omega=f/\tau\left(f\right)$. $\tau^{5}=\mbox{id}$, so $\omega^{5}=1$.\qed

Any polynomial that is $\tau$-semi-invariant and contains a monomial $m$ must contain all the monomials from the $\tau$-orbit of $m$. It is easy to check that the $m_{d,i}$ above are representatives of all orbits of monomials of degree $d\leq4$ in $5$ variables, the result follows.\qed

Now look at how the elements of $D$ act on these polynomials. Since $D$ preserves the basis of $\complex^{5}$, all the monomials are $D$-semi-invariant, so every $\tau$--invariant polynomial must be preserved. Applying $g=\left[p,a_{1},\ldots,a_{5}\right]$ ($p$ prime, $0\leq a_{i}<p$, $a_{i}$ not all equal), get:

\begin{lemma}\label{lemma:sl5:elts_of_D}
 For any $g=\left[n,a_{1},\ldots,a_{5}\right]\in D$, $a_{i}$ not all equal (i.e.\ $g$ is not scalar), the following  hold (replacing $g$ by its scalar multiple if necessary) for some parameter $a\in\integers$ ($0\leq a\leq n$):
\[
\begin{array}{l|l}
A_{1}=0&\\
\hline
B_{1}=0\mbox{ or }n=2&B_{2}=B_{3}=0\\
\hline
C_{1}=0\mbox{ or }n=3&C_{2}=0\mbox{ or }g=\left[11,a,4^{3}a,4^{6}a,4^{9}a,4^{12}a\right]\\
C_{3}=0\mbox{ or }g=\left[11,a,4^{4}a,4^{8}a,4^{12}a,4^{16}a\right]&C_{4}=0\mbox{ or }g=\left[11,a,4^{1}a,4^{2}a,4^{3}a,4^{4}a\right]\\
C_{5}=0\mbox{ or }g=\left[11,a,4^{2}a,4^{4}a,4^{6}a,4^{8}a\right]&C_{6}=C_{7}=0\\
\hline
D_{1}=0\mbox{ or }n\in\defsetshort{2,4}&D_{2}=0\mbox{ or }g=\left[61,a,34^{2}a,34^{4}a,34^{6}a,34^{8}a\right]\\
D_{3}=0\mbox{ or }g=\left[61,a,34^{1}a,34^{2}a,34^{3}a,34^{4}a\right]&D_{4}=0\mbox{ or }g=\left[61,a,34^{4}a,34^{8}a,34^{12}a,34^{16}a\right]\\
D_{5}=0\mbox{ or }g=\left[61,a,34^{3}a,34^{6}a,34^{9}a,34^{12}a\right]&D_{6}=D_{7}=0\mbox{ or }n=2\\ 
D_{8}=0\mbox{ or }g=\left[11,a,4^{1}a,4^{2}a,4^{3}a,4^{4}a\right]&D_{9}=0\mbox{ or }g=\left[11,a,4^{2}a,4^{4}a,4^{6}a,4^{8}a\right]\\
D_{10}=D_{11}=0&D_{12}=0\mbox{ or }g=\left[11,a,4^{3}a,4^{6}a,4^{9}a,4^{12}a\right]\\
D_{13}=0\mbox{ or }g=\left[11,a,4^{4}a,4^{8}a,4^{12}a,4^{16}a\right]&D_{14}=0
\end{array}
\]
\end{lemma}
\pf
This relies on fairly straightforward algebra and using that $\sum_{i}{a_{i}}=0\mod{n}$. All these calculations are almost identical, so only one of them (for $D_{2}\neq0$) will be shown here.

If $D_{2}\neq0$, then the semi-invariance suggests:
\[3a_{1}+a_{2}\equiv3a_{2}+a_{3}\equiv3a_{3}+a_{4}\equiv3a_{4}+a_{5}\equiv3a_{5}+a_{1}\mod n\]
This immediately says $n\neq3$ (otherwise get $a_{1}\equiv\ldots\equiv a_{5}\mod n$, making $g$ a scalar), and hence, by Lemma~\ref{lemma:sl5:makeprime}, $n$ is not divisible by $3$.
Furthermore, it is easy to see that
\[
3a_{1}\equiv 2a_{2}+a_{3},\ 
3a_{2}\equiv 2a_{3}+a_{4},\ 
3a_{3}\equiv 2a_{4}+a_{5},\ 
3a_{4}\equiv 2a_{5}+a_{1},\ 
3a_{5}\equiv 2a_{1}+a_{2} \mod{n} 
\]
Since $a_{1}+\cdots+a_{5}\equiv0\mod{n}$, get 
\begin{eqnarray*}
 0&\equiv&2\left(a_{1}+\cdots+a_{5}\right)\equiv2a_{1}+\left(2a_{2}+a_{3}\right)+a_{3}+\left(2a_{4}+a_{5}\right)+a_{5}\mod{n}\\
&\equiv&2a_{1}+3a_{1}+a_{3}+3a_{3}+a_{5}\equiv5a_{1}+4a_{3}+a_{5}\mod{n}\\
&\equiv&5a_{1}+4a_{3}+\left(3a_{3}-2a_{4}\right)\equiv5a_{1}+7a_{3}-2\left(3a_{2}-2a_{3}\right)\mod{n}\\
&\equiv&5a_{1}+11a_{3}-3\left(2a_{2}\right)\equiv5a_{1}+11a_{3}-3\left(3a_{1}-a_{3}\right)\equiv14a_{3}-4a_{1}\mod{n}
\end{eqnarray*}
giving $4a_{1}\equiv 14a_{3}\mod{n}$. Similarly, get:
\[
4a_{1}\equiv14a_{3},\ 
4a_{2}\equiv14a_{4},\ 
4a_{3}\equiv14a_{5},\ 
4a_{4}\equiv14a_{1},\ 
4a_{5}\equiv14a_{2} 
\]


Since $n$ is not a multiple of $3$, $3$ is invertible $\mod{n}$, and so, writing
\begin{eqnarray*}
 9a_{1}&\equiv&2\left(3a_{2}\right)+3a_{3}\equiv 7a_{3}+2a_{4}\mod{n}\\
 27a_{1}&\equiv&20a_{4}+7a_{5}\mod{n}\\
 81a_{1}&\equiv&61a_{5}+20a_{1}\mod{n}
\end{eqnarray*}
one deduces that either $61|n$ or $a_{1}\equiv a_{5}\mod{n}$. By symmetry (or repeating the calculation for $a_{2},\ldots,a_{5}$) one sees that
\[61a_{1}=61a_{2}=61a_{3}=61a_{4}=61a_{5}\mod{n}\]
and since $n,a_{1},\ldots,a_{5}$ are assumed not to all have a common divisor, one sees that either $n=61$ or $g$ to be a scalar. Since $14\equiv34\cdot4\mod{61}$, the result follows.
\qed

\begin{cor}
 Let $G\subset\slgroup{5}{\complex}$ be a finite irreducible monomial group that induces a non-weakly-exceptional singularity. Then either \abs{D} or $\abs{D}/5$ is in \defsetshort{2^k,3^k,11^k,61^k} for some positive integer $k$.
\end{cor}
Now one needs to look at the possible isomorphism classes of $T$. The remainder of this section will complete the proof of the main technical theorem by excluding most of the possibilities for $T$. In particular, Corollary~\ref{cor:sl5:mono:11-61} will deal with the case where the size of $D$ is divisible by $11$ or $61$, and Proposition~\ref{prop:sl5:nastypolynomial} will show that the remaining groups only need to be checked against the diagonal hypersurfaces.

\begin{cor}
Let $G\subset\slgroup{5}{\complex}$ be a finite irreducible monomial group that induces a non-weakly-exceptional singularity, and $\exists g\in G$ an element of order $11$ or $61$. Then $G=D\rtimes\cyclgroup{5}$ (with $D$ as above).
\end{cor}
\pf It is easy to see that $D\rtimes\cyclgroup{5}\subseteq G$. Assume the inequality is strict. Then by looking at the action of $G$ on the polynomials, it is clear that $C_{i},C_{j}\neq0$ for some $2\leq i\neq j\leq5$. Then any elements of $D$ must be of the form specified in Lemma~\ref{lemma:sl5:elts_of_D}. However, it is easy to see that an element being in to of the forms at the same time means (in the notation of Lemma~\ref{lemma:sl5:elts_of_D}) that $a=0$, and so this is the identity element, leading to a contradiction. A similar argument works for the relevant $D_{i}$-s.\qed

\begin{cor}\label{cor:sl5:mono:11-61}
 If $G$ contains an element of order $11$ or $61$ but induces a singularity that is not weakly-exceptional, then $G$ belongs to one of $16$ conjugacy classes given in Theorem~\ref{thm:sl5:main}(\ref{thm:sl5:main:mono}) (defined by the choice of a primitive root of unity modulo $11$ or $61$ resp.\ and by whether or not $G$ contains non-trivial scalars).
\end{cor}
In Corollary~\ref{cor:sl5:mono:11-61}, the groups with elements of order $11$ are automorphisms of the well-known Klein cubic threefold (see~\cite{Adler78}). Similarly, the groups with elements of order $61$ are automorphisms of the Klein quartic threefold (see~\cite[\S4.3]{Gonzalez11}).

\begin{prop}\label{prop:sl5:nastypolynomial}
 Let $G$ be a finite monomial group as described above preserving the polynomial
\begin{eqnarray*}
h\left(x_{1},\ldots,x_{5}\right)&=&D_{6}\left(x_{1}^{2}x_{2}^{2}+\omega x_{2}^{2}x_{3}^{2}+\omega^{2} x_{3}^{2}x_{4}^{2}+\omega^{3} x_{4}^{2}x_{5}^{2}+\omega^{4} x_{5}^{2}x_{1}^{2}\right)\\&+&D_{7}\left(x_{1}^{2}x_{3}^{2}+\omega x_{2}^{2}x_{4}^{2}+\omega^{2} x_{3}^{2}x_{5}^{2}+\omega^{3} x_{4}^{2}x_{1}^{2}+\omega^{4} x_{5}^{2}x_{2}^{2}\right)
\end{eqnarray*}
semi-invariant for some values of $D_{6},D_{7}$ not both zero, and some $\omega$ a fifth root of $1$. Then $\omega=1$, and the polynomial 
$f\left(x_{1},\ldots,x_{5}\right)=x_{1}^{2}+x_{2}^{2}+x_{3}^{2}+x_{4}^{2}+x_{5}^{2}$
is also $G$-semi-invariant.
\end{prop}
\pf Decompose $G=D\rtimes T$, $\tau\in T$ as above. Lemma~\ref{lemma:sl5:elts_of_D} implies that for any $g\in D$, $g^{2}$ is a scalar, and so such $D$ also leaves $f$ semi-invariant. Therefore, it remains to check that the representatives of generators of $T$ leave $f$ semi-invariant. This is obviously true if $T\cong\cyclgroup{5}$ (then $T$ is generated by the image of $\tau$).

Therefore, it remains to show the proposition holds for $\cyclgroup{5}\subsetneq T\subseteq\symgroup{5}$. Looking at the subgroups of \symgroup{5}, this means $\dihgroup{10}\subseteq T\subseteq\symgroup{5}$. In particular $\exists\delta\in G\setminus D$, such that the image of $\delta$ is (up to conjugation and choosing $\tau$ appropriately) $(2\ 5)(3\ 4)\in\dihgroup{10}\subseteq T\subseteq\symgroup{5}$
Therefore, $\exists\lambda_{i}\in\complex\setminus0$ such that $g$ is defined by \defmap{\left(x_{1},x_{2},x_{3},x_{4},x_{5}\right)}{\left(\lambda_{1}x_{1},\lambda_{5}x_{5},\lambda_{4}x_{4},\lambda_{3}x_{3},\lambda_{2}x_{2}\right)}.

Applying this to $h$ and solving the resulting equations, get
$\lambda_{2}^{2}=\lambda_{1}^{2}\omega^{4}$, $\lambda_{3}^{2}=\lambda_{1}^{2}\omega^{3}$, $\lambda_{4}^{2}=\lambda_{1}^{2}\omega^{3}$, $\lambda_{5}^{2}=\lambda_{1}^{2}\omega$.
By the definition of the semi-direct product, have $\delta^{2}\in D$, and so
$\lambda_{1}^{2}=C\left(-1\right)^{a_{1}},\ \lambda_{3}\lambda_{4}=C\left(-1\right)^{a_{3}}$.
This and the fact that (by construction) $\omega^{5}=1$ implies that $\omega=1$, and hence
$\lambda_{1}^{2}=\lambda_{2}^{2}=\lambda_{3}^{2}=\lambda_{4}^{2}=\lambda_{5}^{2}$,
making $f$ semi-invariant under the action of $\delta$.

Hence the proposition holds unless $\dihgroup{10}\subsetneq T\subseteq\symgroup{5}$. Doing the same calculation (simplified, as $\omega=1$) for the elements of $G\setminus D$ that are preimages of $(1\ 2\ 3)\in\altgroup{5}\subset\symgroup{5}$ and $(2\ 3\ 5\ 4)\in\mathbb{GA}\left(1,5\right)\subset\symgroup{5}$ excludes the remaining $3$ possibilities for $T$.\qed

This concludes the proof of the main technical result of this note, showing that the groups whose size is not divisible by $11$ or $61$ need only be checked against the diagonal hypersurfaces. The final result will estimate the maximal size and the number of groups preserving such hypersurfaces, thus completing the proof of the note's main result.

\begin{prop}
 Assume $G$ is a monomial group leaving the polynomial \[f\left(x_{1},\ldots,x_{5}\right)=x_{1}^{d}+\omega x_{2}^{d}+\omega^{2}x_{3}^{d}+\omega^{3}x_{4}^{d}+\omega^{4}x_{5}^{d}\]
(for some $d\leq4$) semi-invariant. Then $G$ belongs to one of finitely many conjugacy classes.
\end{prop}
\pf Since $\forall g_{D}\in D$, $g_{D}^{5d}=\mbox{id}$ (as $g_{D}^{d}$ is a scalar), there are only finitely many possibilities for $D$ up to choice of basis (in fact, at most $5d^{4}$). The element $\tau$ has been chosen explicitly, so one only needs to worry about elements of $G$ not generated by $D$ and $\tau$. But since $\cyclgroup{5}\subseteq T\subseteq\symgroup{5}$, by Remark~\ref{rmk:sl5:subgroups_of_S5}, $G$ must be generated by $D$, $\tau$ and one more element $\delta$, with the projection of $\delta$ into $\symgroup{5}$ being one of the $4$ known elements. Say $\delta\left(f\right)=\psi f$.

Since $\delta\in\symgroup{5}$, $\delta^{k}=\mbox{id}$ for some $k\leq6$, and so $\psi^{k}=1$. Furthermore, $\forall i\leq5$, $\exists j,l\leq5$ such that $\delta\left(x_{i}^{d}\right)=\psi\omega^{l}x_{j}^{d}$ (as $f$ is preserved), so any non-zero entries in the matrix of $\delta$ must be roots of $1$ of degree at most $k\cdot5\cdot d\leq30d\leq 120$.

Therefore, there are only finitely many possible conjugacy classes for $G$.\qed

\def\cprime{$'$}

\end{document}